\documentclass{amsart}

\usepackage[lite]{amsrefs}

\usepackage{verbatim}
\usepackage{amssymb}
\usepackage{upref}
\usepackage[all]{xy}

\newtheorem{thm}{Theorem}[section]

\newtheorem{prop}[thm]{Proposition}
\newtheorem{obs}[thm]{Observation}

\theoremstyle{definition}
\newtheorem{defn}[thm]{Definition}
\newtheorem{ex}[thm]{Example}
\newtheorem{q}[thm]{Question}

\newtheorem*{notn}{Notation and Terminology}

\newtheorem*{warn}{Warning}

\theoremstyle{remark}

\newcommand{\secref}[1]{Section~\textup{\ref{#1}}}
\newcommand{\thmref}[1]{Theorem~\textup{\ref{#1}}}

\newcommand{\propref}[1]{Proposition~\textup{\ref{#1}}}

\newcommand{\obsref}[1]{Observation~\textup{\ref{#1}}}
\newcommand{\qref}[1]{Question~\textup{\ref{#1}}}

\newcommand{\righttext}[1]{\qquad\text{#1 }}
\renewcommand{\for}{\righttext{for}}

\newcommand{\midtext}[1]{\quad\text{ #1 }\quad}
\renewcommand{\and}{\midtext{and}}

\DeclareMathOperator{\Hom}{Hom}

\newcommand{\N}{\mathbb N}
\newcommand{\Z}{\mathbb Z}

\newcommand{\C}{\mathbb C}

\renewcommand{\H}{\mathcal H}
\newcommand{\G}{\mathcal G}

\renewcommand{\C}{\mathcal C}

\newcommand{\PP}{\mathcal P}

\renewcommand{\bar}{\overline}

\newcommand{\LL}{\Lambda}
\newcommand{\OO}{\Omega}
\renewcommand{\S}{\Sigma}
\newcommand{\DD}{\Delta}

\renewcommand{\a}{\alpha}
\renewcommand{\b}{\beta}

\newcommand{\g}{\gamma}
\renewcommand{\d}{\delta}

\newcommand{\m}{\mu}

\newcommand{\id}{\text{id}}

\renewcommand{\:}{\colon}

\newcommand{\inv}{^{-1}}

\newcommand{\pause}{\renewcommand{\qed}{}\end{proof}}
\newcommand{\resume}[1]{\begin{proof}[Back to the proof of #1]}
\newcommand{\proofof}[1]{\begin{proof}[Proof of #1]}

\newcommand{\smsquare}[4]{
\xymatrix
@C-5pt
@R-5pt
{
.
\ar[d]_#1
&.
\ar[l]_#2
\ar[d]^#3
\\
.
&.
\ar[l]^#4
}
}

\newcommand{\fillsquare}[5]{
\xymatrix
@C-5pt
@R-5pt
{
.
\ar@{}[dr]|#5
\ar[d]_#1
&.
\ar[l]_#2
\ar[d]^#3
\\
.
&.
\ar[l]^#4
}
}

\newcommand{\ff}[1]{\mathcal P(#1)}
\renewcommand{\gg}[1]{\mathcal G(#1)}
\newcommand{\gf}[1]{\pi (#1)}

\begin{document}

\title{Fundamental groupoids of $k$-graphs}

\author[Pask]{David Pask}
\address{School of Mathematical and Physical Sciences
\\University of Newcastle
\\
NSW 2308
\\
Australia}
\email{davidp@maths.newcastle.edu.au}

\author[Quigg]{John Quigg}
\address{Department of Mathematics and Statistics
\\Arizona State University
\\Tempe, Arizona 85287
\\USA}
\email{quigg@math.asu.edu}

\author[Raeburn]{Iain Raeburn}
\address{School of Mathematical and Physical Sciences
\\University of Newcastle
\\
NSW 2308
\\
Australia}
\email{iain@maths.newcastle.edu.au}

\thanks{This research was supported by grants
from the Australian Research Council and the University of Newcastle}




\begin{abstract}
$k$-graphs are higher-rank analogues of directed graphs which were first
developed to provide combinatorial models for operator algebras of
Cuntz-Krieger type.
Here we develop a theory of the fundamental groupoid of a $k$-graph,
and relate it to the fundamental groupoid of an associated graph
called the $1$-skeleton.
We also explore the
failure, in general, of $k$-graphs to faithfully embed into
their fundamental groupoids.
\end{abstract}

\maketitle

\section{Introduction}
\label{sec:Intro}

\emph{$k$-graphs} are combinatorial structures which are
$k$-dimensional analogues of (directed) graphs. They were introduced by
Kumjian and the first author
\cite{kp:kgraph} to help understand work of Robertson and
Steger on higher-rank analogues of the Cuntz-Krieger algebras
\cite{rob-steg}.
Strictly speaking, $k$-graphs are generalizations of path categories: modulo
conventions as regards composition, the $1$-graphs are precisely the path
categories of ordinary graphs.
As the theory of $k$-graphs has developed, the depth of analogy with 
graphs has been remarkable --- it seems that almost every aspect of 
graphs has a valid and interesting analogue for the more general $k$-graphs.

One of the most useful invariants of a graph is its fundamental
groupoid, which
classifies the coverings of the graph, and thereby
gives a
purely combinatorial approach to covering space theory. In a subsequent paper
\cite{pqr:cover} we will develop a theory of coverings for $k$-graphs; in
preparation for this we develop in this paper an elementary theory of
a \emph{fundamental groupoid} $\gg\LL$ of a $k$-graph $\LL$.

One problem with our fundamental groupoid is that there is no good analogue of
the
usual
reduced-word description of elements of $\gg\LL$. Indeed, it is not even
true that the $k$-graph itself embeds in $\gg\LL$. Our main objective,
therefore, is to realize our fundamental groupoid as a well-controlled
quotient of the fundamental groupoid of an associated graph called
the \emph{$1$-skeleton} of $\LL$. This allows us to perform calculations in a
situation where reduced-word arguments are available.

Our approach to realizing the groupoid of a $k$-graph as a quotient of the
groupoid of its 1-skeleton is to first realize the $k$-graph itself as a
quotient of the associated 1-graph. To see how this could be possible, recall
that a $k$-graph $\LL$ is a category with a ``degree functor''
$d\:\LL\to\N^k$ satisfying a certain factorization property (see
\secref{sec:prelim} for the precise definition). The elements of $\LL$ whose
degree is a standard basis vector can be regarded as the edges of the
1-skeleton, and the various factorizations of an arbitrary element of $\LL$
into edge-paths give the desired equivalence relation on the associated
1-graph.

Because every small category is isomorphic to a quotient of a path category,
it will be clear from the proofs that all our results carry over to arbitrary
small categories; however,
we eschew such a generalization since we have no useful applications.

\medskip

After we completed this paper, we learned of the existence of
\cites{bridson, kum:fundgroupoid}.
\cite{bridson}*{Appendix} develops the elementary theory of
the fundamental group of a small category and proves results similar to some
of ours. Bridson and Haefliger concentrate on the fundamental \emph{group} ---
indeed, they stop just short of defining the fundamental groupoid.
\cite{kum:fundgroupoid} develops, in the specific
context of $k$-graphs, the fundamental groupoid and the existence of the
universal covering.
We also wish to thank Kumjian for bringing \cite{bridson} to our attention.

\medskip

We begin in \secref{sec:prelim} by stating our conventions for $k$-graphs and
also for quotients of path categories of graphs.
In \secref{sec:fundamental} we introduce our notion of fundamental groupoids
of $k$-graphs, using a universal construction from category theory, namely
categories of fractions.
In \secref{sec:present k-graph}, we characterize $k$-graphs as quotients of
the path categories of their $1$-skeletons,
and in \secref{sec:present groupoid} we
perform the same job for the associated fundamental groupoids.

In \secref{sec:geometric}, we give a geometric interpretation of our
fundamental groupoid $\gg\LL$, showing that that at least the fundamental
group agrees with
the usual one for a topological space constructed from $\LL$.
\cite{bridson} mentions that in general the fundamental group of a small
category is isomorphic to that of the classifying space, but in the case of
$k$-graphs our geometric realization seems more elementary.

In the last two sections we consider the question of when the canonical
functor $i\:\LL\to\gg\LL$ is injective, and how we might get around the
problem. In \secref{sec:embed}, we give several very simple examples which
illustrate
that
$i$ is often non-injective.
In
\secref{sec:lambda bar}, we discuss the possibility of replacing $\LL$ by its
image $i(\LL)$, in which we actually have $i(\LL)\to\gg{i(\LL)}=\gg\LL$
injective. Unfortunately we have been unable to prove that $i(\LL)$ has the
unique factorization property.

\section{Preliminaries}
\label{sec:prelim}

\subsection*{Small categories}

Let $\C$ be a small category.
We regard $\C$ more as an algebraic structure rather than a set of
objects and morphisms.
For this we identify the objects with certain
idempotent elements of $\C$.
Among our conventions, we write:

\begin{itemize}
\item
$\C^0$ for the set of objects in $\C$;

\item
$s(a)$ for the domain of $a \in \C$, and call it the \emph{source};

\item
$r(a)$ for the codomain of $a \in \C$, and call it the \emph{range};


\item
composition as juxtaposition: $ab=a\circ b$ whenever
$s(a)=r(b)$;


\item
$u\C=\{a\in\C\mid r(a)=u\}$,
$\C v=\{a\in\C\mid s(a)=v\}$, and $u\C v=\Hom(v,u)$
(for $u,v\in\C^0$).
\end{itemize}

In fact, in general we often write composition of maps as juxtaposition,
especially when we are chasing around commutative diagrams.

\subsection*{Path categories}

In order to deal effectively with quotients of 1-graphs, we regard these
quotients as presentations of categories with generators and relations.
More precisely, a 1-graph may be regarded as the free category generated by
the edges and vertices, and the relations will be any set of ordered pairs
of paths generating the equivalence relation defining the quotient.

For a large part of our development
we mainly
adopt the conventions and results of
Schubert's book \cite{schubert}, particularly regarding the use of 
graphs in (small) category theory --- see Chapter 6 of the Schubert book.
Alternative sources are \cites{maclane, rbrown, higgins, lyndon} (particularly
Chapter III of the latter).
However, we occasionally prefer
notation and terminology from the graph algebra
literature (see \cite{kpr}, for example).

A \emph{graph} is a pair $E=(E^0,E^1)$ of sets equipped with two maps
$s,r\:E^1\to E^0$; $E^0$ comprises the \emph{vertices} and $E^1$ the
\emph{edges}, and the \emph{source} and \emph{range} of an edge $e$ are $s(e)$
and $r(e)$, respectively.
A \emph{diagram of type $E$} in a category $\C$ is a graph morphism
from $E$
to the underlying graph of $\C$
(obtained by forgetting the composition);
for simplicity we'll regard the diagram as a map from
$E$ to $\C$.
$E$ embeds a category $\PP(E)$ having the universal property that
for
every diagram $D$ of type $E$ in a category $\C$ there exists a
unique functor $T_D$ making the diagram
\[
\xymatrix
@C-5pt
@R-5pt
{
E \ar@{^(->}[r] \ar[dr]_D
&\ff E \ar[d]^{T_D}
\\
&\C
}
\]
commute.
$\PP(E)$ is the \emph{path category} of $E$,
and the embedding $E\hookrightarrow \PP(E)$ is the \emph{canonical
diagram of type $E$}.

\warn{To be consistent with the $k$-graph literature (e.g., \cite{kp:kgraph})
and Schubert \cite{schubert}, we write
paths of edges in $E$ in the order appropriate for
composition, i.e., if $e_1,\dots,e_n\in E^1$ with $s(e_i)=r(e_{i+1})$ for
$i=1,\dots,n-1$, then $e_1\cdots e_n$ is a path in $\PP(E)$; this is the
opposite of the convention in much of the literature of graphs (e.g.,
\cite{gt}) and graph (operator) algebras (e.g., \cite{kpr}).}

\smallskip

$E\mapsto\ff E$ is functorial from
graphs to small categories. In fact, it is a left adjoint for
the underlying-graph functor.
A \emph{relation}
for $E$ is a pair $(\a,\b)$ of paths in
$\PP(E)$, where $s(\a)=s(\b)$ and $r(\a)=r(\b)$.
We say
a diagram
$D\:E\to\C$ \emph{satisfies} a relation $(\a,\b)$
if $D(\a)=D(\b)$.

\begin{ex}
A diagram of type
$E=\raisebox{1em}{\xymatrix@1@!0{*{.}\ar[d]_a&*{.}\ar[l]_b\ar[dl]^c\\*{.}}}$
in a category $\C$ has the form
\[
\xymatrix
@C-5pt
@R-5pt
{
X
\ar[d]_f
&Y
\ar[l]_g
\ar[dl]^h
\\
Z
}
\]
and satisfies
the relation
$(ab,c)$ if and only if it commutes in the usual sense.
\end{ex}

If $K$ is a set of relations for $E$, we say a diagram of type $E$
\emph{satisfies $K$}
if it satisfies all
relations in $K$.
The quotient of $\PP(E)$ by the smallest equivalence relation containing $K$
is a category which
we denote 
by $\ff E/K$ (by self-explanatory abuse of notation) and call the
\emph{relative path category of $E$ with relations $K$}.
The
quotient map $\ff E\to\ff E/K$ is injective on objects, so we can
identify the objects of $\ff E/K$ with the vertices of $E$.
We call the composition
\[
\DD\:E\hookrightarrow \ff E\to\ff E/K
\]
the \emph{canonical diagram of type $E$ satisfying $K$}; it has the
universal property that
for every diagram
$D\:E\to\C$ satisfying $K$
there exists a unique functor $T_D$
making the diagram
\[
\xymatrix
@C-5pt
@R-5pt
{
E
\ar[r]^-\DD
\ar[dr]_D
&\ff E/K
\ar[d]^{T_D}
\\
&\C
}
\]
commute.

\begin{ex}
If
$E=\raisebox{1em}{\xymatrix@1@!0{*{.}\ar[d]_a&*{.}\ar[l]_b\ar[dl]^c\\*{.}}}$,
$K=\{(ab,c)\}$,
and
$F=\raisebox{1em}{\xymatrix@1@!0{*{.}\ar[d]_a&*{.}\ar[l]_b\\*{.}}}$,
then $\ff E/K\cong\ff F$.
\end{ex}

\subsection*{$k$-graphs}

We adopt the conventions of \cite{kp:kgraph}:
a \emph{$k$-graph} is a category $\LL$
equipped with a \emph{degree functor} $d\:\LL\to\N^k$ satisfying the
following
\emph{factorization property} (also known as \emph{unique
factorization}): for all $\a\in\LL$ and $n,l\in\N^k$ such that
$d(\a)=n+l$ there exist unique $\b,\g\in\LL$ such that $d(\b)=n$,
$d(\g)=l$, and $\a=\b\g$.

If $\LL$ is a $k$-graph, $n\in\N^k$, and $u,v\in\LL^0$, we write
$\LL^n=d\inv(n)$,
$u\LL^n=u\LL\cap\LL^n$, and
$\LL^nv=\LL v\cap\LL^n$.
A \emph{morphism} of $k$-graphs is a degree-preserving functor.

\section{Fundamental groupoids}
\label{sec:fundamental}

Let $\LL$ be a $k$-graph. As a category, $\LL$ has very few invertible
elements: just the vertices.
We shall associate a groupoid to $\LL$, doing as
little damage to $\LL$ as possible, and making all the elements of $\LL$
invertible.

A standard construction in category theory (see
\cite{schubert}*{Section 19.1}, for example) takes any subset $\S$ of
$\LL$ and produces a
new category $\LL[\S\inv]$, called a
\emph{category of fractions},
and a functor $i\:\LL\to\LL[\S\inv]$
such that $i(a)$ is invertible for all $a\in\S$, and having
the universal
property that
if $T:\Lambda\to \C$ is a functor with
$T(a)$ invertible
for all $a\in\S$ then there is a unique functor
$T'$ making the diagram
\[
\xymatrix
@C-5pt
@R-5pt
{
{\LL}
\ar[r]^-{i}
\ar[dr]_{T}
&{\LL[\S\inv]}
\ar[d]^{T'}
\\
&{\C}
}
\]
commute.
Since $\LL$ is small, we can, and shall, take $\S=\LL$,
and we write
\[
\gg\LL=\LL[\LL\inv].
\]
Since $i(\a)$ is invertible for all $\a\in\LL$, there is a
subcategory of $\gg\LL$ which is a groupoid and contains $i(\LL)$; by
universality this subcategory in fact coincides with $\gg\LL$.
Thus $\gg\LL$ is a groupoid, and is generated as a groupoid by $i(\LL)$.
Moreover,
the universal property may be rephrased as follows: for any functor $T$
from $\LL$ into a groupoid $\H$ there exists a unique groupoid morphism
$T'$ making the diagram
\[
\xymatrix
@C-5pt
@R-5pt
{
{\LL}
\ar[r]^-{i}
\ar[dr]_{T}
&{\gg\LL}
\ar[d]^{T'}
\\
&{\H}
}
\]
commute.
The pair $(\gg\LL,i)$ is unique up to isomorphism.

\begin{defn}
With the above notation, $\gg\LL$ is the \emph{fundamental groupoid} of
$\LL$, and $i$ is the \emph{canonical functor}.
\end{defn}

$\LL$ and $\gg\LL$ have the same vertices (objects), and the restriction
$i|\LL^0$ is the identity map. Since $\gg\LL$ is a groupoid, its
hom-set $v\gg\LL v$
at any vertex $v\in\LL^0$ is a group.

\begin{defn}
The \emph{fundamental group} of $\LL$ at a vertex $v$ is $\pi(\LL,v):=v\gg\LL
v$.
\end{defn}

By construction, for any morphism $T\:\LL\to\OO$ between $k$-graphs,
there is a
unique morphism $T_*$ making the diagram
\[
\xymatrix
{
{\LL}
\ar[r]^-{i}
\ar[d]_{T}
&{\gg\LL}
\ar[d]^{T_*}
\\
{\OO}
\ar[r]_-{i}
&{\gg\OO}
}
\]
commute,
because $iT\:\LL\to\gg\OO$ is a functor from $\LL$ to a groupoid.
It follows quickly from this that
$\LL\mapsto\gg\LL$
is functorial
from
$k$-graphs to groupoids.
In fact, this functor extends readily to small categories, and then
it becomes
a left adjoint for the inclusion functor from
groupoids to small categories.

Since $\LL$ generates $\G(\LL)$ as a groupoid, it is easy to see,
for example, that
$\Lambda$ is connected if and only if its fundamental groupoid $\G(\LL)$
is (a small category $\C$ is \emph{connected} if the equivalence relation on
$\C^0$ generated by $\{(u,v)|u\C v\ne\emptyset\}$ is $\C^0\times\C^0$).

\section{Presentation of $k$-graphs}
\label{sec:present k-graph}

In this section we shall give a presentation of any $k$-graph $\LL$ as a
relative path category.

\begin{notn}
We need an unambiguous notation for the standard basis vectors in $\N^k$ ---
the usual $e_i$ is not good for us because it conflicts with typical notation
for edges of a graph --- we use $n_1,\dots,n_k$ for the standard
basis vectors.
\end{notn}

\begin{defn}
The \emph{$1$-skeleton} of $\LL$ is the graph $E$ with $E^0=\LL^0$
and $E^1=\bigcup_{i=1}^k\LL^{n_i}$, and range and source maps inherited from
$\LL$.
\end{defn}

\begin{defn}
If $e$ and $f$ are composable edges in $E$ with orthogonal degrees, the
factorization property of the degree functor gives unique edges $g$ and $h$
such that $ef$ and $gh$ are both \emph{edge-path factorizations} of the same
element of $\LL$, with the degrees interchanged: $d(e)=d(h)$ and $d(f)=d(g)$.
The ordered pair $(ef,gh)$ is a relation on $E$ which we call a
\emph{commuting square}. We let $S$ denote the set of all commuting squares of
$E$.
\end{defn}

We typically visualize
a commuting square $(ef,gh)$ is as
\[
\xymatrix
@C-5pt
@R-5pt
{
.
\ar[d]_{g}
&.
\ar[l]_{h}
\ar[d]^{f}
\\
.
&.
\ar[l]^{e}
}
\]

\begin{prop}
\label{present kgraph}
Let $\LL$ be a $k$-graph, with 1-skeleton $E$ and commuting squares $S$. Then
there is a unique isomorphism $T$ making the diagram
\[
\xymatrix
@C-5pt
@R-5pt
{
E \ar[r]^-\DD \ar@{^(->}[dr]
&\ff E/S \ar[d]^T_\cong
\\
&\LL
}
\]
commute.
\end{prop}

\begin{proof}
The inclusion $E\hookrightarrow \LL$ is a diagram satisfying $S$, so there is
a unique \emph{functor} $T$ making the diagram commute.
In fact, $T$ is the unique factorization through the relative path category
$\ff E/S$ of the unique functor
$R\:\ff E\to\LL$ extending
the diagram $E\to\LL$.
Thus $T$ is surjective.

For $\lambda,\m\in\ff E$, we have $R(\lambda)=R(\m)$ if and only if $\lambda$
and $\m$ are
edge-path factorizations of the same element of $\LL$.
This gives an equivalence relation on $\ff E$.
It follows from the
factorization property of the degree functor and a routine induction argument
that $S$ generates this equivalence relation, and the result
follows.
\end{proof}

\section{Presentation of fundamental groupoids}
\label{sec:present groupoid}

Let $\LL$ be a $k$-graph. In this section we shall give a presentation of the
fundamental groupoid $\gg\LL$ as a relative path category.
The value of this is that
while
$\gg\LL$ is created specifically to have universal
properties,
we have more tools to effectively compute with, and prove
things about, relative path categories.

Schubert would construct the fundamental groupoid
of a small category such as $\LL$
as a relative path category
of a certain augmented version of the underlying
graph of $\LL$,
obtained by adjoining
inverse edges with appropriate relations.
In our case,
\propref{present kgraph} already gives a presentation of
the $k$-graph $\LL$
as $\ff E/S$; it
seems natural to want to work directly with the 1-skeleton
$E$.
In
\propref{present groupoid} below, we obtain the fundamental groupoid of a
$k$-graph $\LL$
as a relative path category of
an augmented version of the 1-skeleton $E$, namely:

\begin{defn}
For each edge $e\in E^1$ introduce a new edge $e\inv$ with source
and range interchanged: $s(e\inv)=r(e)$ and $r(e\inv)=s(e)$, and put
$E\inv=\{e\inv\mid e\in E^1\}$.
Next let $E^+=E\cup E\inv$.  More precisely, the graph $E^+$
has edges $E^1\cup E\inv$ and vertices $E^0$.
We call $E^+$ the \emph{augmented graph} of $E$, and $E\inv$ the
\emph{inverse edges}.
\end{defn}

It is notationally convenient to write
$(e^{-1})^{-1}$ to mean $e$, because then we can write $e^{-1}$
for any $e\in E^1\cup E^{-1}$.

\begin{defn}
Let $C$ denote the set of relations for $E^+$ of the
form $(e\inv e,s(e))$ for $e\in E^1\cup E\inv$. We call $C$ the
set of \emph{cancellation relations for $E$}.
\end{defn}

\begin{obs}
\label{diagram extend}
Every diagram $D$ of type $E$ in a groupoid $\H$ extends uniquely to a
diagram $D^+$ of type $E^+$ in $\H$ satisfying $C$: just put
\[
D^+(e\inv)=D(e)\inv
\for
e\in E^1.
\]
\end{obs}

We also need to know that the canonical functor $i\:\LL\to\gg\LL$ gives rise
to yet another universal property:

\begin{obs}
Let $\LL$ be a $k$-graph with 1-skeleton $E$ and canonical functor
$i\:\LL\to\gg\LL$.
Then for any diagram $D$ of type $E$ in a groupoid $\H$
there exists a unique morphism $T$ making the diagram
\[
\xymatrix
@C-5pt
@R-5pt
{
E \ar[r]^-{i|E} \ar[dr]_D
&\gg\LL \ar[d]^T
\\
&\H
}
\]
commute.
\end{obs}

\begin{thm}
\label{present groupoid}
Let $\LL$ be a $k$-graph, with 1-skeleton $E$, commuting squares $S$,
augmented graph $E^+$, and cancellation relations $C$.
Let $\DD\:E^+\to\ff
{E^+}/(C\cup S)$ be the canonical diagram of type $E^+$ satisfying $C\cup S$.
Then there is a unique
isomorphism $T$ making the diagram
\[
\xymatrix{
E \ar[r]^{i|E} \ar@{^(->}[d]
&\gg\LL \ar[d]^T_\cong
\\
E^+ \ar[r]_-\DD
&\ff {E^+}/(C\cup S)
}
\]
commute.
\end{thm}

\begin{proof}
Note that
by universality
there is a unique \emph{morphism} $T$ making the diagram commute,
since
$\DD|E$ is a diagram of type $E$, satisfying $S$,
in the groupoid
$\ff {E^+}/(C\cup S)$; we must show $T$ is an isomorphism.
Because
$T$ is bijective on units, and
$\DD(E)$
generates $\ff {E^+}/(C\cup S)$ as a groupoid,
$T$ is surjective.

It suffices to show that there exists a
morphism $R$ making the diagram
\[
\xymatrix{
E \ar[r]^{i|E} \ar@{^(->}[d]
&\gg\LL
\\
E^+ \ar[r]_-\DD
&\ff {E^+}/(C\cup S) \ar[u]_R
}
\]
commute,
for if we have such an $R$ then
\[
RTi|E
=R\DD|E
=i|E,
\]
so $RT=\id$ by universality, hence $T$ is injective.

Since $i|E$ is a diagram of type $E$ in a groupoid,
by \obsref{diagram extend}
it extends
uniquely to a diagram $D$ of type $E^+$ satisfying $C$,
so that the diagram
\[
\xymatrix
@C-5pt
@R-5pt
{
E \ar[r]^-{i|E} \ar@{^(->}[d]
&\gg\LL
\\
E^+ \ar[ur]_D
}
\]
commutes.
Because $i|E$
satisfies
$S$, the diagram $D$ satisfies
$C\cup S$.
Thus there exists a unique morphism $R$ making the diagram
\[
\xymatrix{
&\gg\LL
\\
E^+ \ar[ur]^D \ar[r]_-\DD
&\ff {E^+}/(C\cup S) \ar[u]_R
}
\]
commute.
Combining diagrams, because both triangles of the diagram
\[
\xymatrix{
E \ar[r]^-{i|E} \ar@{^(->}[d]
&\gg\LL
\\
E^+ \ar[ur]^D \ar[r]_-\DD
&\ff {E^+}/(C\cup S) \ar[u]_R
}
\]
commute, so does the outer rectangle, as desired.
\end{proof}

Let's look a little more closely at the case $k=1$:
there are no commuting squares, and a 1-graph $\LL$ is not only isomorphic to,
but in fact coincides with, the path category $\ff E$ of its 1-skeleton $E$.
The above theorem gives us an isomorphism
\[
\gg\LL\cong\ff {E^+}/C.
\]
The left-hand side is by definition the fundamental groupoid of the path
category $\ff E$, and is constructed without direct reference to $E$ itself.
On the other hand, the right-hand side is what we expect the fundamental
groupoid of a graph to be: words in edges and their inverses, modulo
cancellation of subwords like $ee\inv$ or $e\inv e$. This latter deserves to
be associated formally with the graph:

\begin{defn}
Let $E$ be a graph, with augmented graph $E^+$ and cancellation
relations $C$. We define the \emph{fundamental groupoid of $E$} to be
\[
\gf E:=\ff {E^+}/C.
\]
\end{defn}

$E\mapsto\gf E$ is
a left adjoint for the underlying graph functor on groupoids.
We can use this to give another interpretation of \thmref{present groupoid}:
let $\LL$ be an arbitrary $k$-graph, with 1-skeleton $E$, commuting squares
$S$, augmented graph $E^+$, and cancellation relations $C$. Relative path
categories can be formed iteratively (``in stages''):
\[
\gg\LL\cong\ff {E^+}/(C\cup S)
\cong\bigl(\ff {E^+}/C\bigr)/S
=\gf E/S,
\]
where we continue our abuse of notation by using ``$/S$'' to
mean we take the
quotient by the equivalence relation generated by $S$.
Note that the commuting squares $S$ pass unaffected into the relative path
category $\ff {E^+}/C$.

\section{Geometric interpretation}
\label{sec:geometric}

The fundamental groupoid $\mathcal G(\Lambda)$ of a $k$-graph $\Lambda$ is a
purely combinatorial object --- even
in the case $k=1$ it is \emph{not}
isomorphic to the (classical) fundamental groupoid of the topological
realization of the graph, because it does not have enough units.
However, it is almost obvious that the (combinatorial) fundamental groupoid of
a graph is
isomorphic to
the \emph{reduction} of the topological fundamental
groupoid to the subset of the unit space consisting of the vertices of the
graph.
If the graph is connected, then the fundamental \emph{groups} of the $1$-graph
and of the topological realization \emph{are} isomorphic.
Something like this persists for $k>1$.
We only sketch the construction --- it's primarily folklore from
algebraic topology.
For simplicity, assume $\Lambda$ is connected.
Let $E$ be the 1-skeleton of $\LL$, with commuting squares $S$.
Then
let $|E|$ be the usual geometric
realization as a topological space, namely as a 1-cell complex where a 1-cell
is attached for each edge of $E$.
Since the $k$-graph $\Lambda$ is connected, so is the topological space $|E|$.

Next, for each commuting square $(\a,\b)\in S$,
attach a 2-cell to $|E|$ along the path $\a\b\inv$ (using the obvious
identification of elements of $\gf E$ with certain continuous paths in $|E|$).
Doing this for all commuting squares,
we get a 2-cell complex $X$.
Now fix a vertex $x\in \Lambda^0$.
It is a
standard fact from algebraic topology that the fundamental group
$\pi_1(X,x)$ of $X$ at the point $x$
is
isomorphic to the quotient of the fundamental group $\pi_1(|E|,x)$
by the normal subgroup generated by loops of the form
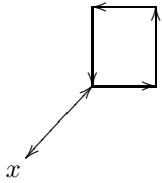
\[
\xymatrix{
& *{} \ar[d]
& *{} \ar[l]
\\
& *{} \ar[r] \ar@{<->}[dl]
& *{} \ar[u]
\\
x
}
\]
i.e., follow any path from $x$ to a corner of one of the attached 2-cells,
then go around the boundary of this 2-cell, and finally ``retrace your steps''
back to $x$.
The boundary of the 2-cell is of the form $\a\b\inv$ for some commuting square
$(\a,\b)\in S$.
It follows from \thmref{present groupoid}, then, that
this quotient group
is also isomorphic to the (combinatorial) fundamental group $\pi(\LL,x)$.
Thus $\pi(\LL,x)$ is isomorphic to the fundamental group of the 2-cell
complex $X$.
Also, the fundamental groupoid $\G(\LL)$ is isomorphic to an appropriate
reduction of the fundamental groupoid of $X$.

We could proceed to attach a 3-cell for each ``commuting cube'' in $\LL$, and
so on up to dimension $k$, but this would have no effect on the fundamental
group of the cell complex.

\section{Failure of embedding}
\label{sec:embed}

Since the canonical functor
$i\:\LL\to\gg\LL$
\emph{is} injective for 1-graphs (folklore),
and $k$-graphs are higher-dimensional generalizations,
it's natural to
ask: Is
$i$ injective for every $k$-graph $\LL$?
It seems that
the standard injectivity result in the
category literature involves a hypothesis which Schubert
\cite{schubert} would call ``$\LL$ admits a calculus of left
fractions'' (some other writers would say ``$\LL$ is localizing''),
which requires that for any $\a,\b\in\LL$ with $s(\a)=s(\b)$ there exist
$\g,\d\in\LL$ such that $\g\a=\d\b$.
This is no help for $\LL$,
since a $k$-graph does not
typically admit a calculus of left fractions.
In fact,
that's a good thing, because the answer to
the injectivity question for $k$-graphs is no!
We shall give counterexamples in this section.
In fact, our
counterexamples will show that injectivity fails as soon as $k\ge 2$.

\begin{ex}
\label{basic}
The following $2$-graph
is the
simplest
example
we could find
where
the canonical functor
$i\:\LL\to\gg\LL$ is noninjective.
It has
a single vertex, 2 edges of degree $(1,0)$, and 3 edges
of degree $(0,1)$, with commuting squares

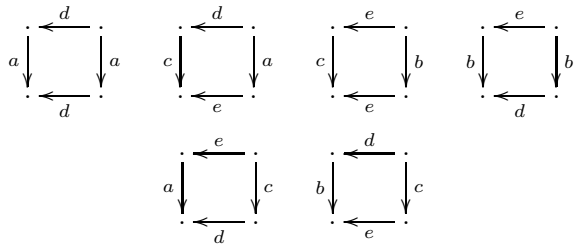
\[
\begin{matrix}
\smsquare{a}{d}{a}{d}
&
\smsquare{c}{d}{a}{e}
&
\smsquare{c}{e}{b}{e}
&
\smsquare{b}{e}{b}{d}
\\
&
\smsquare{a}{e}{c}{d}
&
\smsquare{b}{d}{c}{e}
\end{matrix}
\]
Then in the fundamental group $G$ we have:
\[
a=dad\inv=de\inv c=dbe\inv=b.
\]
In fact, it follows that in $G$ we have
$a=b=c$
and
$d=e$.

Another interpretation of this example: we have a group $G$ with generators
$a,b,c,d,e$ and relations
\[
ad=da,
\quad cd=ea,
\quad ce=eb,
\quad be=db,
\quad ae=dc,
\quad bd=ec,
\]
and the conclusion is that the generators $b,c,e$,
and the last two relations,
are redundant and the group
$G$ is isomorphic to $\Z^2$ (with free commuting generators $a,d$).
\end{ex}

\begin{ex}
In the above example,
we only used the top four commuting squares,
and we tried to
minimize the
number of edges, which
required us to use only 1 vertex.
Here is an example with 4 vertices:
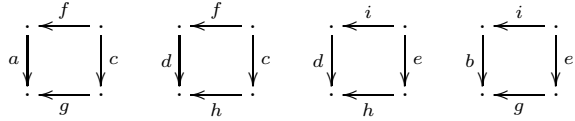
\[
\begin{matrix}
\smsquare{a}{f}{c}{g}
&
\smsquare{d}{f}{c}{h}
&
\smsquare{d}{i}{e}{h}
&
\smsquare{b}{i}{e}{g}
\end{matrix}
\]
where again $i(a)=i(b)$ in the fundamental groupoid.
To embed these
diagrams in a $2$-graph, we have to add more diagrams to ensure that they
give a bijection between the
vertical-horizontal edge-paths and the
horizontal-vertical ones
(see \cite{kp:kgraph}*{Section 6}).
\end{ex}

We could further generalize the previous examples by replacing the edges
$a,b,\dots$ by arbitrary paths $\a,\b,\dots$. Then the squares are interpreted
as commuting diagrams in the $k$-graph --- and for emphasis
we point out
that there is
no reason why $k$ must be 2, except that when
$k>2$ the existence of a $k$-graph
with the desired properties becomes somewhat more delicate:
by \cite{fowlersims_artin}*{Theorem~2.1 and Remark~2.3}
we must also
check consistency of the commuting cubes (``associativity'').

From a $k$-graph $\LL$ we can form $k$ ``component'' 1-graphs
\[
\LL_i:=d\inv(\N n_i)
\qquad,i=1,\dots,k,
\]
where $n_i$ is the $i$th standard basis vector.
The above examples
show
that in general even the component 1-graphs do not embed faithfully in the
$k$-graph fundamental groupoid $\gg\LL$.

Each component 1-graph $\LL_i$ \emph{does} embed faithfully in \emph{its}
fundamental groupoid $\G(\LL_i)$. Thus, we see that in general the $k$-graph
fundamental groupoid $\G(\LL)$ is a kind of ``twisted product'' of nonfaithful
copies of the component 1-graph fundamental groupoids $\G(\LL_i)$.

\section{Lambda bar}
\label{sec:lambda bar}

Let $\LL$ be a $k$-graph,
$d\:\LL\to\N^k$ the degree functor,
and $i:\LL\to \gg \LL$
the canonical functor into the
fundamental groupoid. The universal property of $i$ tells us that there is a
functor $d': \gg \LL \to \Z^k$ such that $d'\circ i=d$. Thus the restriction
$\bar d$ of $d'$ to the image $\bar \LL:=i(\LL)$ is a functor with values in
$\N^k$. We would very much like to know:

\begin{q}
\label{bar kgraph}
Is $(\bar\LL,\bar d)$ a $k$-graph?
\end{q}

The point is that the canonical functor $\bar i$ of $\bar\LL$ into $\G(\bar
\LL)=\G(\LL)$ is injective, so that if the answer is affirmative, we
could for many purposes replace $\LL$ by a $k$-graph which has the same
fundamental groupoid and has $i$ injective.
Fortunately, our general theory does not depend upon the answer to \qref{bar
kgraph}.
Nevertheless, we explore this question a little further:
How could we search for a counterexample?
It suffices to find $\LL$
containing elements $\alpha$ and $\beta$ and factorizations
$\alpha=\gamma\delta$ and $\beta=\epsilon\zeta$ with
\[
i(\alpha)=i(\beta),
\quad
d(\gamma)=d(\epsilon),
\and
i(\gamma)\ne i(\epsilon).
\]
With an eye to keeping things as simple as possible,
suppose
$\gamma$ and $\epsilon$ are edges with degree $(1,0,0)$,
and
$\delta$ and $\zeta$ have degree $(0,1,0)$,
so that
$\alpha$ and $\beta$ have degree $(1,1,0)$
(we'll need the third coordinate shortly).
The techniques we used to find the examples in the preceding section
show that if $\LL$ has commutative diagrams of the form
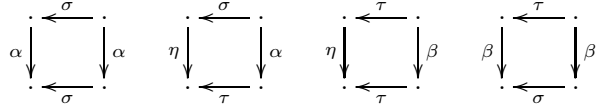
\[
\begin{matrix}
\smsquare{\alpha}{\sigma}{\alpha}{\sigma}
&
\smsquare{\eta}{\sigma}{\alpha}{\tau}
&
\smsquare{\eta}{\tau}{\beta}{\tau}
&
\smsquare{\beta}{\tau}{\beta}{\sigma}
&
\end{matrix}
\]
then $i(\alpha)=i(\beta)$.
Still striving for simplicity, let's consider the possibility that $\sigma$
and $\tau$ are edges, say with degree $(0,0,1)$
(a little thought reveals that if $\sigma$ is an edge with degree the
same as either $d(\gamma)$ or $d(\delta)$
then $\alpha=\beta$, which is not what we want).
If we
label the commuting
squares for $\alpha$, $\beta$, and $\eta$ as
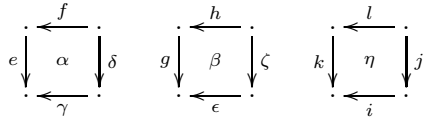
\[
\begin{matrix}
\fillsquare{e}{f}{\delta}{\gamma}{\alpha}
&
\fillsquare{g}{h}{\zeta}{\epsilon}{\beta}
&
\fillsquare{k}{l}{j}{i}{\eta}
\end{matrix}
\]
then the 4 relations give commuting cubes
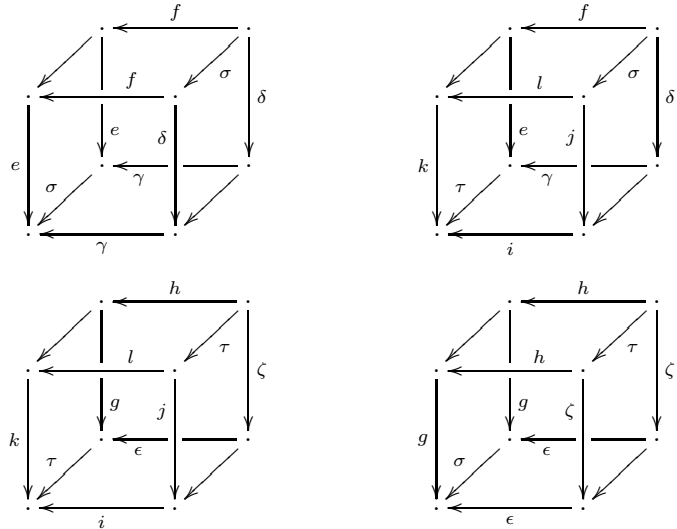
\begin{align*}
\xymatrix
@C-5pt
@R-5pt
{
&.
\ar[dl]
\ar'[d][dd]^e
&&.
\ar[ll]_f
\ar[dl]^\sigma
\ar[dd]^\delta
\\
.
\ar[dd]_e
&&.
\ar[ll]_(.3)f
\ar[dd]_(.3)\delta
\\
&.
\ar[dl]_\sigma
&&.
\ar'[l][ll]^\gamma
\ar[dl]
\\
.
&&.
\ar[ll]^\gamma
}
&&
\xymatrix
@C-5pt
@R-5pt
{
&.
\ar[dl]
\ar'[d][dd]^e
&&.
\ar[ll]_f
\ar[dl]^\sigma
\ar[dd]^\delta
\\
.
\ar[dd]_k
&&.
\ar[ll]_(.3)l
\ar[dd]_(.3)j
\\
&.
\ar[dl]_\tau
&&.
\ar'[l][ll]^\gamma
\ar[dl]
\\
.
&&.
\ar[ll]^i
}
\\
\xymatrix
@C-5pt
@R-5pt
{
&.
\ar[dl]
\ar'[d][dd]^g
&&.
\ar[ll]_h
\ar[dl]^\tau
\ar[dd]^\zeta
\\
.
\ar[dd]_k
&&.
\ar[ll]_(.3)l
\ar[dd]_(.3)j
\\
&.
\ar[dl]_\tau
&&.
\ar'[l][ll]^\epsilon
\ar[dl]
\\
.
&&.
\ar[ll]^i
}
&&
\xymatrix
@C-5pt
@R-5pt
{
&.
\ar[dl]
\ar'[d][dd]^g
&&.
\ar[ll]_h
\ar[dl]^\tau
\ar[dd]^\zeta
\\
.
\ar[dd]_g
&&.
\ar[ll]_(.3)h
\ar[dd]_(.3)\zeta
\\
&.
\ar[dl]_\sigma
&&.
\ar'[l][ll]^\epsilon
\ar[dl]
\\
.
&&.
\ar[ll]^\epsilon
}
\end{align*}
To continue this search for a counterexample
would require answers to the following questions,
which we have so far been unable to supply:

\begin{enumerate}
\item
Can the above cubes be completed to a 3-graph?

\item
Is it true that the relations do not imply $i(\gamma)=i(\epsilon)$?
\end{enumerate}



\end{document}